\documentclass[a4paper, 12pt]{article}
\usepackage{amsmath}
\usepackage{amssymb}
\usepackage{amsthm}
\usepackage{hyperref}
\usepackage{cite}
 \numberwithin{equation}{section}
\newtheorem{thm}{Theorem}
\newtheorem{lemma}{Lemma}

\newtheorem{rem}{Remark}[section]

\date{}

\title { Estimates for lower bounds of eigenvalues of the poly-Laplacian and quadratic polynomial operator of the Laplacian}
\author{Qing-Ming Cheng \  He-Jun Sun \ Guoxin Wei \ Lingzhong Zeng}
\begin{document}
\maketitle

{\narrower\noindent \small { \textsc{Abstract}.}
In this paper, we investigate the Dirchlet eigenvalue problems of poly-Laplacian with any order and quadratic polynomial operator of the Laplacian.
We give some estimates for lower bounds of the sums of their first $k$ eigenvalues which improve the previous results.

}

\footnotetext{



The first author was supported by a Grant-in-Aid for Scientific
Research from JSPS. The second author and the third author were supported by the National Natural Science Foundation of China (Grant Nos.11001130, 11001087).
}

\section{Introduction}

Let $\Omega$ be a bounded domain in an $n$-dimensional
Euclidean space $\mathbb{R}^{n}$, where $n\geq 2$.
The Dirichlet eigenvalue problem of the poly-Laplacian is described by
\begin{equation}
\begin{cases}(-\Delta)^{l}u = \lambda u,\quad \quad \mbox{on}  \ \Omega,\\
u|_{\partial \Omega}=\frac{\partial u}{\partial\nu}|_{\partial \Omega}= \cdots = \frac{\partial^{l-1} u}{\partial\nu^{l-1}}|_{\partial \Omega} =0,
\end{cases}\end{equation}
where $\Delta$ is the Laplacian and $\nu$ denotes the outward unit normal vector field of $\partial\Omega$.
As we known, this problem has a real and discrete spectrum: $0 < \lambda_1 \leq \lambda_2 \leq \cdots \leq \lambda_k \leq \cdots \rightarrow \infty$, where each eigenvalue repeats with its multiplicity.

When $l=1$, problem (1.1) is called the Dirichlet Laplacian problem or the fixed membrane problem.
The asymptotic behavior of its $k$-th eigenvalue $\lambda_k$ relates to geometric properties of $\Omega$ when $k \rightarrow \infty$.
In fact, the following Weyl's asymptotic formula holds
\begin{equation}
\lambda_k  \sim
 \frac{(2 \pi)^2 }{  (\omega_n V(\Omega))^{ \frac{2}{n}} }  k^{\frac{2}{n}}, \quad \mbox{as} \   k \rightarrow \infty,
\end{equation}
where $\omega_n$ denotes the volume of the unit ball in $\mathbb{R}^n$ and $V(\Omega)$ denotes the volume of $\Omega$.
In 1961, P\'{o}lya \cite{P} proved that
\begin{equation}
\lambda_k  \geq
 \frac{(2 \pi)^2 }{  (\omega_n V(\Omega))^{ \frac{2}{n}} }  k^{\frac{2}{n}}
\end{equation}
holds on tiling domains in $\mathbb{R}^{2}$. His proof also works on tiling domains in $\mathbb{R}^{n}$.
Moreover, he conjectured that (1.3) holds for any bounded domain in $\mathbb{R}^n$.
Berezin \cite{Ber} and Lieb \cite{Lie} made some contributions to the partial solution of this conjecture.
In 1983, Li and Yau \cite{LY} proved the following so-called Li-Yau inequality
\begin{equation}
\frac{1}{k} \sum_{j=1}^k  \lambda_j  \geq
\frac{n}{n+2} \frac{(2 \pi)^2}{  (\omega_n V(\Omega))^{\frac{2}{n}} }  k^{ \frac{2}{n}}.
\end{equation}
In 2000, Laptev and Weidl \cite{LW} pointed out that (1.4) can be derived by the Legendre transform of a
result derived by Berezin \cite{Ber}. Hence, (1.4) is also called the Berezin-Li-Yau inequality.
In 2003, adding an additional positive term to the right-hand side of (1.4),  Melas \cite{M} improved (1.4) to
\begin{equation}
\frac{1}{k}   \sum_{j=1}^k  \lambda_j \geq
\frac{n}{n+2} \frac{(2 \pi)^2}{(\omega_nV(\Omega))^{\frac{2}{n}}   }     k^{ \frac{2}{n}}
+ \frac{1}{24(n+2)} \frac{V(\Omega)}{I(\Omega)},
\end{equation}
where  $I(\Omega) = \underset{a \in \mathbb{R}^n}{\textrm{min}}  \int_\Omega |x-a|^2 dx$ is the moment of inertia of $\Omega$.
Recently, Ilyin \cite{I} obtained the following asymptotic lower bound for eigenvalues of problem (1.1):
\begin{equation}
   \begin{aligned}
\frac{1}{k} \sum_{j=1}^k  \lambda_j
 \geq &  \frac{n}{n+2}  \frac{ (2 \pi)^2 }{  (\omega_n V(\Omega))^{\frac{2}{n}} }    k^{ \frac{2}{n}}
  +  \frac{n}{  48}     \frac{ V(\Omega)  }{  I(\Omega)}   \bigg(1- \varepsilon_n(k) \bigg),
  \end{aligned}
\end{equation}
where $0  \leq \varepsilon_n(k) = O( k^{-\frac{2}{n}}) $ is a infinitesimal of $k^{-\frac{2}{n}}$.
Moreover, he derived some explicit inequalities for the particular cases of $n=2, 3, 4$:
\begin{equation}
   \begin{aligned}
\frac{1}{k} \sum_{j=1}^k  \lambda_j
 \geq &  \frac{n}{n+2}  \frac{ (2 \pi)^2 }{ ( \omega_n V(\Omega))^{\frac{2}{n}} }    k^{ \frac{2}{n}}
  +  \frac{n}{  48}   \beta_n  \frac{ V(\Omega)  }{  I(\Omega)},
  \end{aligned}
\end{equation}
where $\beta_2=\frac{119}{120}$, $\beta_3=0.986$ and $\beta_4=0.983$.

When $l=2$, problem (1.1) is called the clamped plate problem.  Agmon \cite{A} and Pleijel \cite{Pl}  obtained
\begin{equation}\lambda_{k}\sim
\frac{(2\pi)^{4}}{ (\omega_{n}V(\Omega)) ^{\frac{4}{n}}  }     k^{\frac{4}{n}}, \quad \mbox{as} \ k
\rightarrow  +\infty.
\end{equation}
In 1985, Levine and Protter \cite{LP}
proved:
\begin{equation}
\frac{1}{k} \sum^{k}_{j=1}\lambda_{j}\geq\frac{n}{n+4}\frac{(2\pi)^{4}}{(\omega_{n} V(\Omega))^{\frac{4}{n}}}    k^{ \frac{4}{n}}.
\end{equation}
For the special case of $n=2$ , Ilyin \cite{I} proved
\begin{equation}
\frac{1}{k} \sum^{k}_{j=1}\lambda_{j}\geq  \frac{16 \pi^2}{3 (V(\Omega))^2}   k^{2} + \frac{12095 \pi}{ 3 \cdot 12096  I(\Omega)} k.
\end{equation}
In 2011, Cheng and Wei \cite{CW1} strengthened (1.9) to
\begin{equation}
   \begin{aligned}
 \frac{1}{k}  \sum^{k}_{j=1}  \lambda_{j}
\geq&  \frac{n}{n+4}\frac{(2\pi)^{4}}   {( \omega_{n}V(\Omega))^{\frac{4}{n}}}      k^{  \frac{4}{n}}\\
&+  \frac{n}{n+2}  \bigg[ \frac{n+2}{12n(n+4)} -\frac{1}{1152n^2(n+4)} \bigg] \frac{(2\pi)^{2}}{(\omega_{n} V(\Omega))^{\frac{2}{n}}   }   \frac{V(\Omega) }{I(\Omega)}   k^{ \frac{2}{n}}\\
&+  \bigg[\frac{1}{576n(n+4)}  - \frac{1}{27648 n^2(n+2)(n+4)}  \bigg]  \Bigg{(}\frac{V(\Omega)}{I(\Omega)}\Bigg{)}^{2}.
  \end{aligned}
\end{equation}

When $l \geq 3$, Levine and Protter \cite{LP} proved
\begin{equation}
\frac{1}{k} \sum^{k}_{j=1}\lambda_{j} \geq  \frac{n}{n+2l}    \frac{(2\pi)^{2l}}{(\omega_{n}V(\Omega))^{\frac{2l}{n}}}   k^{ \frac{2l}{n}}.
\end{equation}
Recently, adding $l$ terms of lower order of $k^{\frac{2l}{n}}$ to its
right-hand side of (1.12), Cheng, Qi and Wei \cite{CQW2011} derived
\begin{equation}
   \begin{aligned}
\frac{1}{k} \sum^{k}_{j=1}\lambda_{j}
&\geq  \frac{n}{n+2l}  \frac{(2\pi)^{2l}}{(\omega_{n}V(\Omega))^{\frac{2l}{n}}}      k^{ \frac{2l}{n}}\\
&+\frac{n}{(n+2l)}\sum_{p=1}^{l}
\frac{l+1-p}{(24)^{p}n\cdots(n+2p-2)}\frac{(2\pi)^{2(l-p)}}{(\omega_{n} V(\Omega))^{\frac{2(l-p)}{n}}}
 \bigg(\frac{V(\Omega)    }{I(\Omega)} \bigg)^{p}       k^{ \frac{2(l-p)}{n}}.
  \end{aligned}
\end{equation}
When $l=1$, (1.13) becomes (1.5).

In this paper, we obtain the following result for problem (1.1).

\begin{thm}
Let $\Omega $ be a bounded domain in an $n$-dimensional
Euclidean space $\mathbb{R}^{n}$. Denote by $\lambda_j$ the $j$-th eigenvalue of problem (1.1). Then we have
\begin{equation}
 \begin{aligned}
 \frac{1}{k}  \sum_{j=1}^k   \lambda_j
 \geq &  \frac{n }{n+2l} \frac{(2\pi)^{2l}}{(\omega_nV(\Omega))^{\frac{2l}{n}} }    k^{ \frac{2l}{n}} \\
&  + \frac{nl }{48} \frac{(2\pi)^{2l-2}  }{( \omega_n V(\Omega) )^{\frac{2l-2}{n}} } \frac{ V(\Omega)   }{I(\Omega)}  k^{ \frac{2l-2}{n}}
  \bigg(1- \varepsilon_n(k) \bigg),
  \end{aligned}
\end{equation}
where $0  \leq \varepsilon_n(k) = O( k^{-\frac{2}{n}}) $ is a infinitesimal of $k^{-\frac{2}{n}}$.
\end{thm}

\begin{rem}
Taking $l=1$ in (1.14), we obtain (1.6).
Moreover, the second term on the right-hand side of (1.13) is
$$
\frac{l }{ 24 (n+2l) }\frac{(2\pi)^{2l-2}}{(\omega_{n}V(\Omega) )^{\frac{2l-2}{n}}}
 \frac{V(\Omega)    }{I(\Omega) }        k^{ \frac{2l-2}{n}}.
$$
Hence, the second term on the right-hand side of (1.14) is $\frac{n(n+2l)}{2}$ times larger than that of (1.13).
Thus, for large $k$, (1.14) is sharper than (1.13).
\end{rem}

Furthermore, we investigate the following Dirichlet eigenvalue problem of quadratic polynomial operator of the Laplacian:
\begin{equation}
\begin{cases}
\Delta^{2}u  - a \Delta u= \Gamma u,\quad \quad \mbox{on}  \ \Omega,\\
u|_{\partial \Omega}=\frac{\partial u}{\partial\nu}|_{\partial \Omega}=0,
\end{cases}\end{equation}
where $a$ is a nonnegative constant. Levine and Protter \cite{LP} proved that the eigenvalues of this problem satisfy
\begin{equation}
   \begin{aligned}
 \Gamma_k \geq&   \frac{n}{n+4}  \frac{ (2 \pi)^4 }{  (\omega_n  V(\Omega))^{\frac{4}{n}} }     k^{ \frac{4}{n}}
  + \frac{na}{n+2}   \frac{(2 \pi)^2}{ (\omega_nV(\Omega))^{\frac{2}{n}}}    k^{\frac{2}{n}}.
  \end{aligned}
\end{equation}

In this paper, we derive the following results for problem (1.15).

\begin{thm}
Let $\Omega $ be a bounded domain in $\mathbb{R}^{n}$. Denote by $\Gamma_j$  the $j$-th eigenvalue of problem (1.15). Then we have
\begin{equation}
   \begin{aligned}
\frac{1}{k}  \sum_{j=1}^k  \Gamma_j
\geq&   \frac{n}{n+4}  \frac{ (2 \pi)^4 }{ ( \omega_nV(\Omega))^{\frac{4}{n}} }    k^{ \frac{4}{n}}
  + \bigg( \frac{n}{  24}\frac{ V(\Omega) }{  I(\Omega)} + \frac{na}{n+2} \bigg)    \frac{(2 \pi)^2}{( \omega_n V(\Omega))^{\frac{2}{n}}}  k^{ \frac{2}{n}}\\
  & +\bigg[ -\frac{ n(n^2-4)  }{3840 }      \frac{ V(\Omega)  }{  I(\Omega)}
  +  \frac{na}{  48  }  \bigg]  \frac{ V(\Omega)  }{  I(\Omega)}    \bigg(1- \varepsilon_n(k) \bigg),
  \end{aligned}
\end{equation}
where $0  \leq \varepsilon_n(k) = O( k^{-\frac{2}{n}}) $ is a infinitesimal of $k^{-\frac{2}{n}}$.
\end{thm}

For the special cases of $n=2, 3, 4$, we prove the following sharper result:

\begin{thm}
Denote by $\Gamma_j$ the $j$-th eigenvalue of problem (1.15) on a bounded domain  $\Omega$ in $\mathbb{R}^{n}$, where $n=2, 3, 4$. Then we have
\begin{equation}
   \begin{aligned}
\frac{1}{k}  \sum_{j=1}^k  \Gamma_j
 \geq&
\frac{n}{n+4}  \frac{ (2 \pi)^4 }{ ( \omega_nV(\Omega))^{\frac{4}{n}} }     k^{ \frac{4}{n}}
 +\bigg(  \frac{n}{  24} \alpha_n  \frac{ V(\Omega) }{  I(\Omega)}
  +   \frac{na }{n+2}  \bigg) \frac{ (2 \pi)^2 }{ ( \omega_nV(\Omega))^{\frac{2}{n}} }   k^{ \frac{2}{n}} \\
&  +  \frac{na }{  48  }  \beta_n \frac{ V(\Omega)  }{  I(\Omega)},
    \end{aligned}
\end{equation}
where $\alpha_2= \frac{12095}{12096}$, $\beta_2= \frac{119}{120}$, $\alpha_3 = 0.991$, $\beta_3 = 0.986$, $\alpha_4=0.985$ and $\beta_4=0.983$.
\end{thm}

Making a modification in the proof of Theorem 3, we can get the following result:

\begin{thm}
Denote by $\Gamma_j$ the $j$-th eigenvalue of problem (1.15) on a bounded domain  $\Omega$ in $\mathbb{R}^{n}$, where $n= 3, 4$. Then we have
\begin{equation}
   \begin{aligned}
\frac{1}{k}  \sum_{j=1}^k  \Gamma_j
\geq&   \frac{n}{n+4}  \frac{ (2 \pi)^4 }{ ( \omega_nV(\Omega))^{\frac{4}{n}} }    k^{ \frac{4}{n}}
  + \bigg( \frac{n}{  24}\frac{ V(\Omega) }{  I(\Omega)} + \frac{na}{n+2} \bigg)    \frac{(2 \pi)^2}{( \omega_n V(\Omega))^{\frac{2}{n}}}  k^{ \frac{2}{n}}\\
  & +\bigg[ -\frac{ n(n^2-4)  }{3840 }      \frac{ V(\Omega)  }{  I(\Omega)}
  +  \frac{na}{  48  } \beta_n \bigg]  \frac{ V(\Omega)  }{  I(\Omega)}.
  \end{aligned}
\end{equation}
\end{thm}

\begin{rem}
Taking $a=0$ in (1.17), (1.18) and (1.19), we can get some results for the clamped plate problem.
\end{rem}

\section{Proofs of the main results}

In order to prove Theorem 1, we need the following lemma derived by Ilyin \cite{I}.
\begin{lemma}
Let
$$
\Psi_s(r)  =
\left \{\aligned
&M,            &  \mbox{for}  \quad  0  \leq r \leq s; \quad \quad \  \\
&M -  L(r-s),  &  \mbox{for}   \quad  s\leq r \leq s + \frac{M}{L};\\
&0,            &  \mbox{for}  \quad   r \geq s + \frac{M}{L}.  \quad \quad
\endaligned\right.
$$
Suppose that $\int_0^{+\infty}  r^b \Psi_s(r)  dr = m^*$ and $d \geq b$.
Then for any decreasing and absolutely continuous function $F$ satisfying the conditions
\begin{equation}
 0 \leq F \leq M,   \quad  \int_0^{+\infty}   r^b   F(r) dr  =  m^*,   \quad  0 \leq   -F'  \leq L,
 \end{equation}
the following inequality holds:
\begin{equation}
\int_0^{+\infty}   r^d   F(r) dr  \geq   \int_0^{+\infty}   r^d  \Psi_s(r) dr.
\end{equation}
\end{lemma}

Now we give the proof of Theorem 1.

{\bf Proof of Thereom 1}\hspace{0.2cm}
Let $u_j$ be an orthonormal eigenfuction corresponding to the $j$-th eigenvalue $\lambda_j$ of problem (1.1).
Denote by $\widehat{u}_{j}(\xi)$ the Fourier transform of  $u_{j}(x)$, which is defined by
\begin{equation}
\widehat{u}_{j}(\xi)= (2\pi)^{-\frac{n}{2}} \int_{\Omega}  u_{j}(x) e^{ i x \cdot \xi}  dx.
\end{equation}
It follows from Plancherel's Theorem that
\begin{equation}
\int_{\Omega} \widehat{u}_{j}(\xi) \widehat{u}_{q}(\xi) d \xi = \delta_{jq}.
\end{equation}

Set
$ h(\xi) = \sum_{j=1}^k  |\widehat{u}_{j}(\xi)|^2.$
From (2.4) and  Bessel's inequality, one can get
\begin{equation}
  h(\xi) = \sum_{j=1}^k  |\widehat{u}_{j}(\xi)|^2  \leq    (2\pi)^{-n}   \int_{\Omega}  | e^{ i x \cdot \xi}|^2 dx = (2\pi)^{-n}  V(\Omega).
\end{equation}
Moreover,  Parsevel's identity implies that
\begin{equation}
\int_{\mathbb{R}^n}  h(\xi)  d\xi = \sum_{j=1}^k \int_{\Omega}   |u_{j}(x)|^2  dx  = k.
\end{equation}
Since
 $$  \nabla  \widehat{u}_{j}(\xi) =   (2\pi)^{-\frac{n}{2}}   \int_{\Omega}  i x u_j(x) e^{ i x \cdot \xi} dx,$$
we have
\begin{equation}
 \sum_{j=1}^k  |\nabla  \widehat{u}_{j}(\xi)|^2 \leq   (2\pi)^{-n} \int_{\Omega} | i x e^{ i x \cdot \xi}|^2 dx  = (2\pi)^{-n}  I(\Omega).
\end{equation}
It follows from (2.5) and (2.7) that
\begin{equation}
  |\nabla  h(\xi)|
  \leq
      2 \big(  \sum_{j=1}^k  | \widehat{u}_{j}(\xi)|^2  \big )^{\frac{1}{2}} \big(  \sum_{j=1}^k  | \nabla \widehat{u}_{j}(\xi)|^2  \big )^{\frac{1}{2}}
   \leq   2 (2\pi)^{-n} \sqrt{V(\Omega) I(\Omega)}.
\end{equation}
Denote by $h^{*}(\xi)= \psi(|\xi|)$ the symmetric decreasing rearrangement (see \cite{Ban, P2}) of $h$.
From
$$
k  =\sum_{j=1}^k \int_{\Omega}   |u_{j}(x)|^2  dx  =\int_{\mathbb{R}^n}  h(\xi)  d\xi
=  \int_{\mathbb{R}^n}  h^*(\xi)  d\xi
 =  n \omega_n  \int_0^{+\infty}  r^{n-1} \psi(r) dr,
$$
we get
\begin{equation}
 \int_0^{+\infty}  r^{n-1} \psi(r) dr  =   \frac{k}{n \omega_n }.
\end{equation}
At the same time, using integration by parts and Parseval's identity, we have
\begin{equation}
   \begin{aligned}
 \int_{\mathbb{R}^n}   |\xi|^{2l}  h(\xi) d\xi
  =& \sum_{j=1}^k  \sum_{p_1, \cdots, p_l=1}^n \int_{\mathbb{R}^n} \bigg| (2\pi)^{-\frac{n}{2}}
  \int_{\Omega}  \xi_{p_1} \cdots \xi_{p_l}    u_j(x)   e^{  i x \cdot \xi }d x \bigg|^2 d \xi  \\
 =& \sum_{j=1}^k  \sum_{p_1, \cdots, p_l=1}^n \int_{\mathbb{R}^n} \bigg| (2\pi)^{-\frac{n}{2}}
  \int_{\Omega}   \frac{\partial^l u_j(x)}{\partial x_{p_1} \cdots \partial x_{p_l} }  e^{  i x \cdot \xi }d x \bigg|^2 d \xi  \\
 =& \sum_{j=1}^k  \sum_{p_1, \cdots, p_l=1}^n \int_{\mathbb{R}^n} \bigg|\frac{\widehat{\partial^l u_j(\xi)}}{\partial x_{p_1} \cdots \partial x_{p_l} } \bigg |^2 d \xi \\
  = &   \sum_{j=1}^k  \sum_{p_1, \cdots, p_l=1}^n \int_{\mathbb{R}^n} \bigg(\frac{ \partial^l u_j(x) }{\partial x_{p_1} \cdots \partial x_{p_l} }  \bigg)^2 d x \\ =& \sum_{j=1}^k  \int_{\Omega}   u_j(x) (- \Delta)^l u_j(x) dx.
 \end{aligned}
\end{equation}
Thus, it yields
\begin{equation}
\sum_{j=1}^k   \lambda_j  =  \int_{\mathbb{R}^n}  |\xi|^{2l}  h(\xi) d\xi.
\end{equation}
Making use of (2.11) and the properties of symmetric decreasing rearrangement, we obtain
\begin{equation}
   \begin{aligned}
\sum_{j=1}^k   \lambda_j  =&  \int_{\mathbb{R}^n}   |\xi|^{2l}  h(\xi) d\xi \geq \int_{\mathbb{R}^n}   |\xi|^{2l}  h^*(\xi) d\xi
                        = n \omega_n  \int_0^{+\infty}  r^{n+2l-1} \psi(r) dr.
   \end{aligned}
\end{equation}

Noticing (2.5), (2.8) and (2.9), we can apply Lemma 1 to $\psi$ with $b=n-1$ and $d=n+2l-1$.
Therefore, using (2.12), we have
\begin{equation}
   \begin{aligned}
  \sum_{j=1}^k   \lambda_j \geq &n \omega_n  \int_0^{+\infty}   r^{n+2l-1} \psi_(r) dr  \geq  n \omega_n  \int_0^{+\infty}  r^{n+2l-1} \Psi_s(r) dr
   \end{aligned}
\end{equation}
with $M= (2\pi)^{-n}  V(\Omega)$, $m_*=  \frac{k}{n \omega_n }$ and $ L = 2 (2\pi)^{-n} \sqrt{V(\Omega) I(\Omega)}$.
Set $t= \frac{Ls}{M}$.
Combining (2.9) and
$$
 \int_0^{+\infty}  r^{n-1} \psi(r) dr = \int_0^{+\infty}  r^{n-1} \Psi_s(r) dr
 =  \frac{M^{n+1}}{n (n+1) L^{n} }    \bigg[ (t+1)^{n+1}  -t^{n+1} \bigg],
$$
it yields
\begin{equation}
  (t+1)^{n+1}  -t^{n+1}   =  k_*,
\end{equation}
where
$$k_*  =  k  \frac{(n+1) L^n}{ \omega_n  M^{n+1}  }.$$
Set $\eta= t  - \frac{1}{2}$. Then (2.14) becomes
\begin{equation}
 (\eta+\frac{1}{2})^{n+1}  -(\eta-\frac{1}{2})^{n+1}  =  k_*.
\end{equation}
The asymptotic expansion for the unique positive root of (2.15) is
\begin{equation}
\eta( k_*)
=    \zeta - \frac{n-1}{24} \zeta^{-1}+ \frac{ (n-1) (n-3)(2n+1)}{5760} \zeta^{-3} + \cdots,
\end{equation}
where $\zeta=  \big(\frac{k_*}{ n+1} \big)^{\frac{1}{n}}$. Then we can deduce
\begin{equation}
   \begin{aligned}
&  \bigg(t(k^*) + 1\bigg)^{n+2l+1} -t(k^*)^{n+2l+1} \\
 =&   \binom{n+2l+1}{1}    \zeta^{n+2l}
  +  2 \bigg[ \frac{1}{2^3}
\binom{n+2l+1}{3}   -  \frac{n-1}{48} \binom{n+2l+1}{2}  \binom{2}{1}     \bigg] \zeta^{n+2l-2}\\
  &+ 2\bigg[   \frac{1}{2^5}  \binom{n+2l+1}{5}   -\frac{1}{2^3}\frac{n-1}{24}\binom{n+2l+1}{4}  \binom{4}{1}
 + \frac{1}{2} \frac{(n-1)^2}{ 24^2 } \binom{n+2l+1}{3}  \binom{3}{1}  \\
  & +  \frac{1}{2}\frac{(n-1)(n-3)(2n+1)}{5760} \binom{n+2l+1}{2}  \binom{2}{1}      \bigg] \zeta^{n+2l-4}
 +\cdots\\
 =&   (n+2l+1)  \bigg[  \zeta^{n+2l}    +  \frac{l (n+2l) }{12}   \zeta^{n+2l-2}  + \frac{(n+2l)C(n,l)}{5760}  \zeta^{n+2l-4}  +\cdots \bigg],
  \end{aligned}
\end{equation}
where $\binom{q}{t} = \frac{q!}{ t! (q-t)!}$ and
$$
   \begin{aligned}
C(n, l)=&
 (n+2l-1) \bigg[  (n+2l-2)( 6l-7n+1)  +5(n-1)^2   \bigg]\\
 &+(n-1)(n-3)(2n+1).
   \end{aligned}
$$
Using (2.17), we get
\begin{equation}
   \begin{aligned}
& n \omega_n  \int_0^{+\infty}  r^{n+2l-1} \Psi_s(r) dr   \\
=&   \frac{n \omega_n M^{n+2l+1}}{(n+2l) (n+2l+1) L^{n+2l} }    \bigg[ \big(t(k_*)+1 \big)^{n+2l+1}  -t(k_*)^{n+2l+1} \bigg]\\
= &     \frac{n \omega_n M^{n+2l+1}}{(n+2l) L^{n+2l}}
\bigg[   \big(\frac{k_*}{ n+1} \big)^{\frac{n+2l}{n}}
+  \frac{l (n+2l) }{12}   \big(\frac{k_*}{ n+1} \big)^{\frac{n+2l-2}{n}}\\
 & +  \frac{(n+2l)C(n,l )}{5760}  \big(\frac{k_*}{ n+1} \big)^{\frac{n+2l-4}{n}}  +\cdots \bigg].
  \end{aligned}
\end{equation}
Substituting $k_* =  k  \frac{(n+1) L^n}{ \omega_n  M^{n+1}  } $, $M= (2\pi)^{-n}  V(\Omega)$ and $ L = 2 (2\pi)^{-n} \sqrt{V(\Omega) I(\Omega)}$
into (2.18), we have
\begin{equation}
 \begin{aligned}
  & n \omega_n  \int_0^{+\infty}  r^{n+2l-1} \Psi_s(r) dr\\
  = &  \frac{n}{n+2l} \omega_n^{-\frac{2l}{n}}  M^{-\frac{2l}{n}}  k^{1+\frac{2l}{n}}
  + \frac{nl}{12} \omega_n^{-\frac{2l-2}{n}}  \frac{ M^{2-\frac{2l-2}{n}}  }{L^2}  k^{1+\frac{2l-2}{n}}\\
&  + \frac{n C(n, l)}{5760}  \omega_n^{-\frac{2l-4}{n}}   \frac{ M^{4- \frac{2l-4}{n}} }{L^4}  k^{1+\frac{2l-4}{n}} +O( k^{1+\frac{2l-6}{n}} )\\
= &  \frac{n }{n+2l} \frac{(2\pi)^{2l} }{ (\omega_nV(\Omega))^{ \frac{2l}{n}} }    k^{1+\frac{2l}{n}}
  + \frac{nl }{48} \frac{(2\pi)^{2l-2}  }{  (\omega_n  V(\Omega))^{ \frac{2l-2}{n}}  }  \frac{ V(\Omega)  }{I(\Omega)}  k^{1+\frac{2l-2}{n}}\\
  & + \frac{  n C(n, l)}{92160} \frac{(2\pi)^{2l-4} }{  ( \omega_n  V(\Omega))^{ \frac{2l-4}{n}} }  \bigg(\frac{ V(\Omega)  }{I(\Omega) } \bigg)^2  k^{1+\frac{2l-4}{n}}
  + O( k^{1+\frac{2l-6}{n}} ).
  \end{aligned}
\end{equation}
Inserting (2.19) into (2.13), we know that (1.14) is true. This completes the proof of Theorem 1.
$\hfill{} \Box$

\vskip 3mm

{\bf Proof of Thereom 2}\hspace{0.2cm}
It follows from (2.10) that
\begin{equation}
   \begin{aligned}
\sum_{j=1}^k  \Gamma_j  =&   \sum_{j=1}^k   \int_{\Omega}  u_j(x)  \bigg( \Delta^2 u_j(x) - a \Delta u_j(x) \bigg) d x\\
                        = & \int_{\mathbb{R}^n}   |\xi|^{4}  h(\xi) d\xi  +a  \int_{\mathbb{R}^n}   |\xi|^{2}  h(\xi) d\xi\\
                        \geq& \int_{\mathbb{R}^n}   |\xi|^{4}  h^*(\xi) d\xi +a \int_{\mathbb{R}^n}   |\xi|^{2}  h^*(\xi) d\xi \\
                        = &   n \omega_n \bigg(  \int_0^{+\infty}  r^{n+3} \psi(r) dr  + a    \int_0^{+\infty}  r^{n+1} \psi(r) dr \bigg).
    \end{aligned}
\end{equation}
Then, applying Lemma 1 to $ \psi$ and using (2.20), we obtain
\begin{equation}
   \begin{aligned}
\sum_{j=1}^k  \Gamma_j \geq&   n \omega_n  \bigg( \int_0^{+\infty}  r^{n+3} \Psi_s(r) dr  + a   \int_0^{+\infty}  r^{n+1} \Psi_s(r) dr  \bigg).
    \end{aligned}
\end{equation}
Observe that $C(n,l)=-24n^2 +96$ when $l=2$ and $C(n,l)=- 4(3n+2) (n-1)$ when $l=1$.
Therefore, from (2.19),  we have
\begin{equation}
   \begin{aligned}
 & n \omega_n  \bigg( \int_0^{+\infty}  r^{n+3} \Psi_s(r) dr  + a   \int_0^{+\infty}  r^{n+1} \Psi_s(r) dr  \bigg)\\
=&   \frac{n}{n+4}  \frac{ (2 \pi)^4 }{ ( \omega_nV(\Omega))^{\frac{4}{n}} }    k^{1+\frac{4}{n}}
  + \bigg( \frac{n}{  24}\frac{ V(\Omega) }{  I(\Omega)} + \frac{na}{n+2} \bigg)   \frac{(2 \pi)^2}{ (\omega_nV(\Omega))^{\frac{2}{n}}}   k^{1+\frac{2}{n}}\\
  & +\bigg[ -\frac{ n(n^2-4)  }{3840 }      \frac{ V(\Omega)  }{  I(\Omega)} + a \frac{n}{  48  }  \bigg]  \frac{ V(\Omega)  }{  I(\Omega)}  k
   + O(   k^{ 1-\frac{2}{n} } ).
    \end{aligned}
\end{equation}
Then it is easy to find that (1.17) holds. This completes the proof of Theorem 2.
$\hfill{} \Box$

\vskip 3mm
{\bf Proof of Thereom 3}\hspace{0.2cm}
When $n=2$, making use of (1.7) and (1.10), we have
\begin{equation}
   \begin{aligned}
\sum_{j=1}^k  \Gamma_j
 = & \int_{\mathbb{R}^n}   |\xi|^{4}  h(\xi) d\xi  +a  \int_{\mathbb{R}^n}   |\xi|^{2}  h(\xi) d\xi\\
\geq&   2 \omega_2    \int_0^{+\infty}  r^{5} \Psi_s(r) dr  +  2 a \omega_2  \int_0^{+\infty}  r^{3} \Psi_s(r) dr \\
\geq&  \frac{1}{3}  \frac{ (2 \pi)^4 }{  (\omega_2  V(\Omega))^2 }   k^{3}
 +\bigg(  \frac{ \alpha_2}{  12 I(\Omega)}
  +   \frac{  a}{2V(\Omega) }   \bigg) \frac{ (2 \pi)^2 }{  \omega_2  } k^{2}
   +  \frac{ a }{  24 }  \beta_2 \frac{ V(\Omega)  }{  I(\Omega)} k,
    \end{aligned}
\end{equation}
where $\alpha_2= \frac{12095}{12096}$ and $\beta_2= \frac{119}{120}$.

When $n=3$, it follows from (2.21) that
\begin{equation}
   \begin{aligned}
\sum_{j=1}^k  \Gamma_j \geq&   3 \omega_3  \int_0^{+\infty}  r^{6} \Psi_s(r) dr  +  3 a\omega_3  \int_0^{+\infty}  r^{4} \Psi_s(r) dr.
    \end{aligned}
\end{equation}
Now we make an estimate for the lower bound of $\int_0^{+\infty}  r^{6} \Psi_s(r) dr$.
Since
$$
   \begin{aligned}
 \int_0^{+\infty}  r^{6} \Psi_s(r) dr
=   \frac{ M^{8}}{56 L^{7} }    \bigg[ \big(t(k_*)+1 \big)^{8}  -t(k_*)^{8} \bigg],
  \end{aligned}
$$
we need to estimate $\big(t(k_*)+1 \big)^{8}  -t(k_*)^{8} $.
The equation (2.14) becomes $(t+1)^4 -t^4 = k_*$ when $n=3$.
Its positive root $t(k_*)$ is
$$
t(k_*) = \frac{1}{2} \big(\rho(k_*)- \varrho(k_*) \big) -\frac{1}{2},
$$
where $\rho(k_*)= \big( k_* +  \sqrt{k_*^2 +\frac{1}{27}} \big)^{\frac{1}{3}} $
and $\varrho(k_*)= \big(- k_* +  \sqrt{k_*^2 +\frac{1}{27}} \big)^{\frac{1}{3}} $.
Set
$ \vartheta (k_*) =     \frac{1}{2} \big(\rho(k_*)- \varrho(k_*) \big)$.
Then we have
\begin{equation}
   \begin{aligned}
 \big(  t(k_*)+1 \big)^8  -t(k_*)^8
=&  8 \vartheta(k_*)^7 + 14 \vartheta(k_*)^5  +\frac{7}{2}  \vartheta(k_*)^3   +\frac{1}{8}   \vartheta(k_*)\\
=  &  \frac{1}{16} \bigg[\big(\rho(k_*)- \varrho(k_*) \big)^7 +  7 \big(\rho(k_*)- \varrho(k_*) \big)^5   \\
& \quad \ +7 \big(\rho(k_*)- \varrho(k_*)\big)^3 +   \big(\rho(k_*)- \varrho(k_*) \big) \bigg].
  \end{aligned}
\end{equation}
Observe that
\begin{equation}
 \rho(k_*)\cdot \varrho(k_*) =\frac{1}{3}.
 \end{equation}
Then, using (2.26), we have
\begin{equation}
   \begin{aligned}
&\big(\rho(k_*)- \varrho(k_*) \big)^7 \\
=&    \rho(k_*) \big(\rho(k_*)^6 +7  \varrho(k_*)^6 \big) +21 \rho(k_*)^2 \varrho(k_*)^2 \big( \rho(k_*)^3 - \varrho(k_*)^3\big)\\
&
  -35 \rho(k_*)^3 \varrho(k_*)^3 \big( \rho(k_*) - \varrho(k_*) \big)- \varrho(k_*) \big(7 \rho(k_*)^6 +\varrho(k_*)^6 \big)\\
=&   \big( k_* +  \sqrt{k_*^2 +\frac{1}{27}} \big)^{\frac{1}{3}}  \big( 16 k_*^2   -12 k_*  \sqrt{k_*^2 +\frac{1}{27}}   -1 \big) +\frac{14}{3}  k_* \\
   &  - \big(- k_* +  \sqrt{k_*^2 +\frac{1}{27}} \big)^{\frac{1}{3}}  \big( 16 k_*^2  +12 k_*  \sqrt{k_*^2 +\frac{1}{27}} -1\big),
  \end{aligned}
 \end{equation}
\begin{equation}
   \begin{aligned}
& \big(\rho(k_*)- \varrho(k_*) \big)^5\\
=&   \rho(k_*)^5     -5  \rho(k_*) \varrho(k_*) \big( \rho(k_*)^3 - \varrho(k_*)^3 \big)    +10\rho(k_*)^2 \varrho(k_*)^2 \big( \rho(k_*)- \varrho(k_*)  \big) -  \varrho(k_*)^5\\
=&   \big( k_* +  \sqrt{k_*^2 +\frac{1}{27}} \big)^{\frac{5}{3}} -\big(- k_* +  \sqrt{k_*^2 +\frac{1}{27}} \big)^{\frac{5}{3}} -\frac{10}{3} k_*  +\frac{10}{9}  \big( k_* +  \sqrt{k_*^2 +\frac{1}{27}} \big)^{\frac{1}{3}}
\\
&  -\frac{10}{9}  \big( -k_* +  \sqrt{k_*^2 +\frac{1}{27}} \big)^{\frac{1}{3}}
  \end{aligned}
 \end{equation}
 and
\begin{equation}
   \begin{aligned}
&\ 7 \big(\rho(k_*)- \varrho(k_*)\big)^3 +   \big(\rho(k_*)- \varrho(k_*) \big) \\
=  & 7 \big[ \rho(k_*)^3 - \varrho(k_*)^3 -3\rho(k_*)  \varrho(k_*) \big(\rho(k_*) -  \varrho(k_*) \big)  \big] +   \big(\rho(k_*)- \varrho(k_*) \big) \\
= &  14k_*   -6  \big( k_* +  \sqrt{k_*^2 +\frac{1}{27}} \big)^{\frac{1}{3}}  + 6 \big(- k_* +  \sqrt{k_*^2 +\frac{1}{27}} \big)^{\frac{1}{3}}.
  \end{aligned}
 \end{equation}
Substituting (2.27-2.29) into (2.25), we obtain
\begin{equation}
   \begin{aligned}
&\big(  t(k_*)+1 \big)^8  -t(k_*)^8 \\
= &   \big( k_* +  \sqrt{k_*^2 +\frac{1}{27}} \big)^{\frac{1}{3}}  \big(  k_*^2 - \frac{3}{4} k_* \sqrt{k_*^2 +\frac{1}{27}} \big)
 - \frac{7}{16}  \big( -k_* +  \sqrt{k_*^2 +\frac{1}{27}} \big)^{\frac{5}{3}} \\
 &+\bigg[ \frac{ 7}{16} \big( k_* +  \sqrt{k_*^2 +\frac{1}{27}} \big)^{\frac{5}{3}}  -  \big(- k_* +  \sqrt{k_*^2 +\frac{1}{27}} \big)^{\frac{1}{3}}  \big(  k_*^2 + \frac{3}{4} k_* \sqrt{k_*^2 +\frac{1}{27}} \big) \bigg] \\
 & + \frac{7}{144} \bigg[ \big( k_* +  \sqrt{k_*^2 +\frac{1}{27}} \big)^{\frac{1}{3}}  -    \big(- k_* +  \sqrt{k_*^2 +\frac{1}{27}} \big)^{\frac{1}{3}} \bigg]- \frac{7}{24}  k_*.
  \end{aligned}
\end{equation}
Now we make some estimates for some terms in the right hand side of (2.30).
The first term can be estimated as follows:
\begin{equation}
   \begin{aligned}
 \big( k_* +  \sqrt{k_*^2 +\frac{1}{27}} \big)^{\frac{1}{3}}  \big(  k_*^2 - \frac{3}{4} k_* \sqrt{k_*^2 +\frac{1}{27}} \big)
\geq &  \frac{2^{\frac{1}{3}}}{ 4}  k_*^{\frac{7}{3}}  -\frac{2^{\frac{1}{3}}}{72} k_*^{\frac{1}{3}}.
  \end{aligned}
\end{equation}
Here we use the inequality $\sqrt{k_*^2 +\frac{1}{27}} \leq k_* + \frac{1}{54 k_*}$ since $k_*$ is large.
The second term is
\begin{equation}
-\frac{7}{16} \big(- k_* +  \sqrt{k_*^2 +\frac{1}{27}} \big)^{\frac{5}{3}} \geq  - \frac{7 \cdot 2^{\frac{1}{3}} }{16  \cdot 54 \cdot 18} k_*^{-\frac{5}{3}}.
\end{equation}
The third term is
\begin{equation}
   \begin{aligned}
& \frac{7}{16}  \big( k_* +  \sqrt{k_*^2 +\frac{1}{27}} \big)^{\frac{5}{3}}  -\big( -k_* +  \sqrt{k_*^2 +\frac{1}{27}} \big)^{\frac{1}{3}}  \big(  k_*^2 + \frac{3}{4} k_* \sqrt{k_*^2 +\frac{1}{27}} \big)\\
\geq &   \frac{7 \cdot 2^{\frac{2}{3}}}{ 12} k_*^{\frac{5}{3}}    -\frac{2^{\frac{2}{3}}}{432} k_*^{-\frac{1}{3}}.
  \end{aligned}
\end{equation}
The fourth term is
\begin{equation}
   \begin{aligned}
\frac{7}{144} \bigg[  \big( k_* +  \sqrt{k_*^2 +\frac{1}{27}} \big)^{\frac{1}{3}}  - \big( -k_* +  \sqrt{k_*^2 +\frac{1}{27}} \big)^{\frac{1}{3}} \bigg]
 \geq & \frac{7\cdot 2^{\frac{1}{3}} }{144} k_*^{\frac{1}{3}} - \frac{7 \cdot 2^{\frac{2}{3}}  }{16 \cdot 54} k_*^{-\frac{1}{3}}.
  \end{aligned}
\end{equation}
Therefore, using (2.31-2.34) in (2.30), we have
\begin{equation}
   \begin{aligned}
\big(  t(k_*)+1 \big)^8  -t(k_*)^8
\geq &  \frac{2^{\frac{1}{3}}}{ 4}  k_*^{\frac{7}{3}}  +  \frac{7 \cdot 2^{\frac{2}{3}}}{ 12} k_*^{\frac{5}{3}}  - \frac{7}{24}  k_*
 +\frac{ 5 \cdot 2^{\frac{1}{3}}}{144} k_*^{\frac{1}{3}} -\frac{2^{\frac{2}{3}}}{96} k_*^{-\frac{1}{3}}\\
  &- \frac{7 \cdot 2^{\frac{1}{3}} }{16  \cdot 54 \cdot 18} k_*^{-\frac{5}{3}} \\
  \geq &  \frac{2^{\frac{1}{3}}}{ 4}  k_*^{\frac{7}{3}}  +  \frac{7 \cdot 2^{\frac{2}{3}}}{ 12} k_*^{\frac{5}{3}}  - \frac{7}{24}  k_*.
    \end{aligned}
\end{equation}
Here we used the fact that $k_* \geq 1$. In fact, noticing $ k_*  \geq \frac{ (n+1)(4\pi)^n}{\omega_n^2} (\frac{n}{n+2})^{\frac{n}{2}}$, it is not difficult to observe that $k_* = 4 k  L^3 ( \omega_3 )^{-1} M^{-4} \geq  \tau:= \frac{432 \sqrt{15} \pi}{25} \approx 210.25$ when $n=3$. Hence, when
$ \alpha \geq \frac{7}{24}  \tau^{-\frac{2}{3}}$, the following inequality
$$
\alpha k_*^{\frac{5}{3}} \geq  \frac{7}{24}  k_*,
$$
holds for $ k_* \in [\tau, +\infty)$.
Since $1 -   \frac{6}{7}\cdot 2^{\frac{1}{3}}\alpha \leq 1- \frac{1 }{4} \cdot 2^{\frac{1}{3}} \tau^{-\frac{2}{3}}  \approx 0.9911$, we can conclude that
\begin{equation}
   \begin{aligned}
\big(  t(k_*)+1 \big)^8  -t(k_*)^8
\geq &    \frac{2^{\frac{1}{3}}}{ 4}  k_*^{\frac{7}{3}}  +  \frac{7 \cdot 2^{\frac{2}{3}}}{ 12} \alpha_3  k_*^{\frac{5}{3}},
    \end{aligned}
\end{equation}
where
$\alpha_3 = 0.991.$
Therefore, using (2.36), we derive
\begin{equation}
   \begin{aligned}
3 \omega_3  \int_0^{+\infty}  r^{6} \Psi_s(r) dr   = &   \frac{3 \omega_3  M^{8}}{56 L^{7} }    \bigg[ \big(  t(k_*)+1 \big)^8  -t(k_*)^8  \bigg]\\
     \geq &       \frac{3 \cdot 2^{\frac{1}{3}} \omega_3  M^{8}}{224 L^{7} }     k_*^{\frac{7}{3}}
   +     \frac{    2^{\frac{2}{3}}  \omega_3  M^{8}}{32 L^{7} }    \alpha_3  k_*^{\frac{5}{3}}  \\
  =&  \frac{3}{7}  \frac{(2\pi)^4}{( \omega_3V(\Omega))^{\frac{4}{3}}  }     k^{\frac{7}{3}}
  + \frac{ 1   }{8} \alpha_3   \frac{(2\pi)^2}{  (\omega_3 V(\Omega))^{  \frac{2}{3}  } } \frac{  V(\Omega)    }{ I(\Omega)} k^{ \frac{5}{3} }.
    \end{aligned}
\end{equation}
At the same time, it follows from (1.7) that
\begin{equation}
   \begin{aligned}
 3 \omega_3  \int_0^{+\infty}  r^{4} \Psi_s(r) dr
  \geq &   \frac{3}{5}   \frac{(2\pi)^2}{ (\omega_3 V(\Omega))^{ \frac{2}{3}} }  k^{\frac{5}{3}}
   + \frac{ 1   }{16} \beta_3     \frac{ V(\Omega)   }{I(\Omega)} k,
    \end{aligned}
\end{equation}
where $\beta_3 = 0.986.$
Substituting (2.37) and (2.38) into (2.24), we obtain
\begin{equation}
   \begin{aligned}
\sum_{j=1}^k  \Gamma_j \geq&      \frac{3}{7}  \frac{(2\pi)^4}{( \omega_3V(\Omega))^{\frac{4}{3}}  }     k^{\frac{7}{3}}
  + \bigg( \frac{  1    }{8} \alpha_3   \frac{  V(\Omega)   }{ I(\Omega)}
  +  \frac{3a}{5}    \bigg) \frac{(2\pi)^2}{ (\omega_3 V(\Omega) )^{ \frac{2}{3}} } k^{ \frac{5}{3} }\\
&  + \frac{   a    }{16}  \beta_3   \frac{ V(\Omega)   }{I(\Omega)} k.
    \end{aligned}
\end{equation}

When $n=4$, it follows from (2.21) that
\begin{equation}
   \begin{aligned}
\sum_{j=1}^k  \Gamma_j \geq&   4 \omega_4  \int_0^{+\infty}  r^{7} \Psi_s(r) dr  +  4 a\omega_4  \int_0^{+\infty}  r^{5} \Psi_s(r) dr.
    \end{aligned}
\end{equation}
Now we make an estimate for the lower bound of $\int_0^{+\infty}  r^{7} \Psi_s(r) dr$.
Since
$$
   \begin{aligned}
 \int_0^{+\infty}  r^{7} \Psi_s(r) dr
=   \frac{ M^{9}}{72 L^{8} }    \bigg[ \big(t(k_*)+1 \big)^{9}  -t(k_*)^{9} \bigg],
  \end{aligned}
$$
we need to estimate $\big(t(k_*)+1 \big)^{9}  -t(k_*)^{9} $.
The equation (2.14) becomes $(t+1)^5 -t^5 = k_*$ when $n=4$.
Its positive root $t(k_*)$ is
$$
t(k_*) =  \theta(k_*) -\frac{1}{2},
$$
where $\theta(k_*)=    \sqrt{\frac{ \sqrt{20k_* +5}}{10} - \frac{1}{4} } $.
Then we have
\begin{equation}
   \begin{aligned}
 &\big(  t(k_*)+1 \big)^9  -t(k_*)^9\\
=&  9 \theta(k_*)^8 + 21 \theta(k_*)^6  +\frac{63}{8}  \theta(k_*)^4   +\frac{9}{16}   \theta(k_*)^2  +  \frac{1}{2^8}\\
=& \frac{9}{25} k_*^2 + \frac{6}{25} k_* \sqrt{20k_*+5}  -\frac{18}{25} k_*  + \frac{3}{50} \sqrt{20k_*+5} -\frac{7}{50}\\
\geq &  \frac{9}{25} k_*^2 + \frac{12 \sqrt{5}}{25} k_*^{\frac{3}{2}}    -\frac{18}{25} k_* .
  \end{aligned}
\end{equation}
Here we used the fact that $k_* \geq 1$. In fact, noticing $ k_*  \geq \frac{ (n+1)(4\pi)^n}{\omega_n^2} (\frac{n}{n+2})^{\frac{n}{2}}$, it is not difficult to observe that $k_* = 5 k  L^4 ( \omega_4 )^{-1} M^{-5} \geq  \sigma:= \frac{5 \cdot 2^{12} }{9} \approx 2275.56 $ when $n=4$. Hence, when
$ \alpha \geq \frac{18}{25}  \sigma^{-\frac{1}{2}}$, the following inequality
$$
\alpha k_*^{\frac{3}{2}} \geq  \frac{18}{25}  k_*,
$$
holds for $ k_* \in [\sigma, +\infty)$.
Since $1 -   \frac{5  \sqrt{5}}{12}  \alpha \leq 1- \frac{3 \sqrt{5} }{10}   \sigma^{-\frac{1}{2}}  \approx 0.9859$, we can conclude that
\begin{equation}
   \begin{aligned}
\big(  t(k_*)+1 \big)^9  -t(k_*)^9
\geq &    \frac{9}{25} k_*^2 + \frac{12 \sqrt{5}}{25} \alpha_4 k_*^{\frac{3}{2}},
    \end{aligned}
\end{equation}
where
$\alpha_4 = 0.985.$
Therefore, using (2.42), we deduce
\begin{equation}
   \begin{aligned}
4 \omega_4  \int_0^{+\infty}  r^{7} \Psi_s(r) dr   = &   \frac{ \omega_4  M^{9}}{18 L^{8} }    \bigg[ \big(  t(k_*)+1 \big)^9  -t(k_*)^9  \bigg]\\
     \geq &       \frac{  \omega_4  M^{9}}{50 L^{8} }     k_*^{2}
   +     \frac{    2 \sqrt{5}  \omega_4  M^{9}}{75 L^{8} }    \alpha_4  k_*^{\frac{3}{2}}  \\
  =&  \frac{1}{2}  \frac{(2\pi)^4}{  \omega_4V(\Omega)    }     k^{2}
  + \frac{ 1   }{6} \alpha_4   \frac{(2\pi)^2}{  (\omega_4 V(\Omega))^{  \frac{1}{2}  } } \frac{  V(\Omega)    }{ I(\Omega)} k^{ \frac{3}{2} }.
    \end{aligned}
\end{equation}
Meanwhile, from (1.7), we have
\begin{equation}
   \begin{aligned}
 4 \omega_4  \int_0^{+\infty}  r^{5} \Psi_s(r) dr
  \geq &   \frac{2}{3}   \frac{(2\pi)^2}{ (\omega_4 V(\Omega))^{ \frac{1}{2}} }  k^{\frac{3}{2}}
   + \frac{ 1   }{12} \beta_4    \frac{ V(\Omega)   }{I(\Omega)} k,
    \end{aligned}
\end{equation}
where $\beta_4 = 0.983.$
Substituting (2.43) and (2.44) into (2.40), we obtain
\begin{equation}
   \begin{aligned}
\sum_{j=1}^k  \Gamma_j \geq&      \frac{1}{2}  \frac{(2\pi)^4}{  \omega_4V(\Omega)    }     k^{2}
  + \bigg( \frac{ 1   }{6} \alpha_4  \frac{  V(\Omega)    }{ I(\Omega)}
  +  \frac{2a}{3}    \bigg) \frac{(2\pi)^2}{  (\omega_4 V(\Omega))^{  \frac{1}{2}  } }  k^{ \frac{3}{2} }\\
&  + \frac{   a    }{12}  \beta_4   \frac{ V(\Omega)   }{I(\Omega)} k.
    \end{aligned}
\end{equation}

Therefore, synthesizing (2.23), (2.39) and (2.45), we conclude that (1.18) is true. This concludes the proof of Theorem 3.
$\hfill{} \Box$

\vskip 3mm
{\bf Proof of Thereom 4}\hspace{0.2cm}
When $n=3$, using (2.35), we derive
\begin{equation}
   \begin{aligned}
3 \omega_3  \int_0^{+\infty}  r^{6} \Psi_s(r) dr
     \geq &       \frac{3 \cdot 2^{\frac{1}{3}} \omega_3  M^{8}}{224 L^{7} }     k_*^{\frac{7}{3}}
   +     \frac{    2^{\frac{2}{3}}  \omega_3  M^{8}}{32 L^{7} }       k_*^{\frac{5}{3}}
   -\frac{ \omega_3 M^8}{64 L^7} k_* \\
  =&  \frac{3}{7}  \frac{(2\pi)^4}{( \omega_3V(\Omega))^{\frac{4}{3}}  }     k^{\frac{7}{3}}
  + \frac{ 1   }{8}    \frac{(2\pi)^2}{  (\omega_3 V(\Omega))^{  \frac{2}{3}  } } \frac{  V(\Omega)    }{ I(\Omega)} k^{ \frac{5}{3} }
  - \frac{1}{256} \bigg( \frac{V(\Omega)}{I(\Omega)} \bigg)^2 k
  .
    \end{aligned}
\end{equation}
Substituting (2.38) and (2.46) into (2.24), we have
\begin{equation}
   \begin{aligned}
\sum_{j=1}^k  \Gamma_j \geq&
\frac{3}{7}  \frac{(2\pi)^4}{( \omega_3V(\Omega))^{\frac{4}{3}}  }     k^{\frac{7}{3}}
  + \bigg(\frac{ 1   }{8} \frac{  V(\Omega)    }{ I(\Omega)}  + \frac{3a}{5}\bigg)
  \frac{(2\pi)^2}{  (\omega_3 V(\Omega))^{  \frac{2}{3}  } }  k^{ \frac{5}{3} }\\
 & +\bigg(- \frac{1}{256}  \frac{V(\Omega)}{I(\Omega)}
   + \frac{ a  }{16} \beta_3   \bigg )  \frac{ V(\Omega)   }{I(\Omega)} k.
    \end{aligned}
\end{equation}
When $n=4$,
it follows from (2.41) that
\begin{equation}
   \begin{aligned}
4 \omega_4  \int_0^{+\infty}  r^{7} \Psi_s(r) dr
     \geq &       \frac{  \omega_4  M^{9}}{50 L^{8} }     k_*^{2}
   +     \frac{    2 \sqrt{5}  \omega_4  M^{9}}{75 L^{8} }       k_*^{\frac{3}{2}}   - \frac{ \omega_4 M^9}{ 25 L^8} k_*\\
  =&  \frac{1}{2}  \frac{(2\pi)^4}{  \omega_4V(\Omega)    }     k^{2}
  + \frac{ 1   }{6}     \frac{(2\pi)^2}{  (\omega_4 V(\Omega))^{  \frac{1}{2}  } } \frac{  V(\Omega)    }{ I(\Omega)} k^{ \frac{3}{2} }
  - \frac{1}{80} \bigg( \frac{V(\Omega)}{I(\Omega)} \bigg)^2 k.
    \end{aligned}
\end{equation}
Substituting (2.44) and (2.48) into (2.40), we obtain
\begin{equation}
   \begin{aligned}
\sum_{j=1}^k  \Gamma_j
\geq &  \frac{1}{2}  \frac{(2\pi)^4}{  \omega_4V(\Omega)    }     k^{2}
  + \bigg(\frac{ 1   }{6} \frac{  V(\Omega)    }{ I(\Omega)} + \frac{2a}{3} \bigg)  \frac{(2\pi)^2}{  (\omega_4 V(\Omega))^{  \frac{1}{2}  } }  k^{ \frac{3}{2} }\\
 &+\bigg( - \frac{1}{80}   \frac{V(\Omega)}{I(\Omega)}
   + \frac{ a   }{12} \beta_4   \bigg) \frac{ V(\Omega)   }{I(\Omega)} k
  .
    \end{aligned}
\end{equation}
Therefore, combining (2.47) and (2.49), we conclude that (1.19) is true. This completes the proof of Theorem 4.
$\hfill{} \Box$



{\small
\begin{flushleft}
\textsc{Qing-Ming Cheng\\
Department of Mathematics, Graduate School of Science and Engineering, Saga University, Saga 840-8502, Japan}\\
\textit{E-mail address:} \verb"cheng@ms.saga-u.ac.jp"
\end{flushleft}

\begin{flushleft}
\textsc{He-Jun Sun\\
Department of Applied Mathematics, College of Science, Nanjing University of Science and Technology,
 Nanjing 210094, P. R. China}\\
\textit{E-mail address:} \verb"hejunsun@163.com"
\end{flushleft}

\begin{flushleft}
\textsc{Guoxin Wei\\
School of Mathematical Sciences, South China Normal University, Guangzhou 510631, P. R. China}\\
\textit{E-mail address:} \verb"weigx03@mails.tsinghua.edu.cn"
\end{flushleft}

\begin{flushleft}
\textsc{Lingzhong Zeng\\
Department of Mathematics,  Graduate School of Science and
Engineering, Saga University, Saga 840-8502, Japan}\\
\textit{E-mail address:} \verb"lingzhongzeng@yeah.net"
\end{flushleft}

}
\end{document}